\documentclass[11 pt]{article}

\usepackage{amsmath,amssymb,amsthm}
\usepackage{mathrsfs}
\usepackage{amsfonts}
\usepackage{geometry}

\newtheorem{theorem}{Theorem}

\newtheorem{conjecture}[theorem]{Conjecture}
\newtheorem{lemma}[theorem]{Lemma}

\newtheorem{remark}[theorem]{Remark}

\newcommand{\R}{\mathbb{R}}

\newcommand{\N}{\mathbb{N}}
\newcommand{\Z}{\mathbb{Z}}
\newcommand{\seteq}{:=}

\DeclareMathOperator{\Fix}{Fix}
\DeclareMathOperator{\Supp}{Supp}

\newcommand{\rwr}{\wr_r}

\title{An exploration of normalish subgroups of R. Thompson's groups $F$ and $T$}
\author{Collin Bleak}
\date{\today}
\begin{document}
\maketitle
\abstract{In this short note, we show that R. Thompson's group $F$ admits a normalish amenable subgroup, and that the standard copy of $F$ in R. Thompson's group $T$ is normalish in $T$.  We further conjecture that if $F$ is non-amenable, then $T$ does not admit a normalish amenable subgroup, and therefore that the reduced $C^*$ algebra of $T$ is in fact simple in that case.}

\vspace{.1 in}


\section{Introduction}
In \cite{HaagerupOlesen}, the authors show that if R. Thompson's group $T$ has a simple reduced $C^*$ algebra ($T$ is \emph{$C^*$-simple}), then $F$ is non-amenable.  Later, in \cite{BleakJuschenko} it is shown that the Kesten Test (commonly used to detect Powers' Criterion for $C^*$-simplicity) cannot be used in the case of the group $T$.  This and other anecdotal evidence has lead some to speculate that $T$ might fail to be $C^*$-simple despite the fact that $T$ is a group with trace.  If true, $T$ would provide a new and interesting example of this rare phenomenon (see \cite{LeBoudec} for Adrien Le Boudec's examples of such groups).

The paper \cite{BKKO} provides a new test for $C^*$-simplicity for a group.  Namely, if a group admits no \emph{normalish} amenable subgroups, then it is $C^*$-simple.  A group \emph{$M$ is normalish in a group $G$} if and only if given any finite set $K:=\{k_1,k_2,\ldots k_m\}\subset G$ we have

\[
\cap_{i=1}^m M^{k_i} \neq \{1_G\}.
\]

Focussing now on the R. Thompson groups, consider $S^1=\R/\Z$ (we use the parameterisation provided by the map $t\mapsto e^{2\pi i t}$).  The \emph{standard R. Thompson's groups $F<T$} are groups of orientation-preserving piecewise-linear homeomorphisms of $S^1$ which preserve the dyadic rationals and which admit at most finitely many breaks in slope,  appearing only over the dyadic rationals $\left\{\frac{a}{2^k}\mid a,k\in\Z\right\}$, and with all slopes of affine components being integral powers of $2$.  In this case, $T$ is the full group of such homeomorphisms, whereas $F$ is the subgroup of $T$ which stabilises the point $0$ under the natural action of $T$ on $S^1$.  See \cite{CFP} for a survey of the important R. Thompson groups.

Our first result has bearing on R. Thompson's group $F$.
\begin{theorem}\label{FNormalish}
R. Thompson's group $F$ admits a normalish amenable subgroup.
\end{theorem}
Our second observation is that the canonical version of $F$ in $T$ is a normalish subgroup in R. Thompson's group $T$.

\begin{remark}\label{TNormalishIfFAmenable}
For the standard R. Thompson groups $F<T$, we have that $F$ is normalish in $T$.
\end{remark}

In particular, in the case that $F$ is amenable, then $T$ admits a normalish amenable subgroup.  This of course agrees with the Haagerup-Olesen result that if $T$ is $C^*$-simple, then $F$ is non-amenable.

{\flushleft {\it Acknowledgements:}} We would like to thank Emmanuel Breuillard and Kate Juschenko for interesting discussions of the $C^*$-simplicity of the groups $F$ and $T$, where their questions and discussions motivated the creation of this note.

\subsection{Further explorations}
It is of abiding interest whether the converse of the Haagerup-Olesen result is also true.  As mentioned before, some researchers have noted that the result in \cite{BleakJuschenko} provides a first indication that $T$ may not be $C^*$-simple even if $F$ is non-amenable.  We do not take this view.  

The set $\mathscr{X}$ of subgroups of $T$ which admit no embedded copies of R. Thompson's group $F$, and no embedded non-abelian free subgroups, is a very interesting and complex set.  Two key results which may be of assistance in understanding constraints on the elements of groups in the set $\mathscr{X}$ are Brin's Ubiquity Theorem, and Lemma 1.9 of \cite{BKMStructure}.  To state Brin's Ubiquity Theorem, we require Brin's notion of \emph{orbital}.  We say a group $H\leq$ PL$_o$($I$) has an orbital $(a,b)$ if $(a,b)$ is an open interval in $I:=[0,1]$ and the interval $(a,b)$ is a component of support of the action of $H$ on $I$.  Likewise, if $g\in H$, we say an interval $(c,d)$ is an orbital of $g$ if it is an orbital of $\langle g\rangle$.  If $H$ has an orbital $(a,b)$ and $c\in(a,b)$,  and $g\in H$ has an orbital of the form $(a,c)$ or $(c,b)$ then we say \emph{$g$ approaches $a$ (or respectively $b$) in $(a,b)$}. 

\begin{theorem}[Brin's Ubiquity]
Let $H$ be a subgroup of PL$_o$($I$).  Assume that $H$ has an orbital $(a, b)$
and that some element of $H$ approaches one end $a$ or $b$ in $(a, b)$ but not the other.
Then $H$ contains a subgroup isomorphic to $F$.
\end{theorem}

A form of the contrapositive of Lemma 1.9 of \cite{BKMStructure} is the following lemma.
\begin{lemma}Let $1<k\in \N$ and suppose $\mathcal{F}\seteq\left\{g_1,g_2,\ldots, g_k\right\}$ is a set of orientation-preserving homeomorphisms of $S^1$ each of which admits a non-trivial fixed set, and suppose further that 
\[
\cap_{1\leq i\leq k}\Fix(g_i) =\emptyset,
\]
then $\langle g_1,g_2,\ldots,g_k\rangle$ contains non-abelian free subgroups. 
\end{lemma}

Together, these results imply that for any  $H\in \mathscr{X}$ the elements of $H$ admit very strong conditions on how their components of support can overlap.  While we will not investigate this further here, we note that our proof that $F$ is normalish in $T$ does not appear to translate directly to apply to such a group $H$, and so we conjecture the following.
\begin{conjecture}
Every normalish subgroup of $T$ contains embedded copies of R. Thompson's group $F$ or of non-abelian free subgroups.
\end{conjecture}

\section{Normalish subgroups exist}
In this section we show our Theorem  \ref{FNormalish} and Remark \ref{TNormalishIfFAmenable}.  First though, we recall some notation.

Let $g\in F\cup T$, and recall the definition
\[
\Supp(g)\seteq\{x\in (0,1) \mid x\cdot b \neq x\}.
\]

We are now ready to prove Theorem \ref{FNormalish}.
\begin{proof}Consider the standard generator $x_0$ of $F$, and another function $b\in F$ so that $\Supp(b) = (1/4,1/2)$.  Since $1/4\cdot x_0 = 1/2$, we have that $W=\langle x_0,b\rangle\cong \Z\rwr \Z$, noting that this is a standard construction.  We observe that the base group $D$ of $W$ is generated by $X\seteq\{b^{x_0^k}\mid k\in \Z\}$ which is all of the conjugates of $b$ by $x_0$.  Further note that 
\[
D\cong \bigoplus _\Z \Z
\] as the generators in the set $X$ are pairwise disjoint.

We now observe that as any element of $F$ is linear near some neighbourhood of $0$, with slope some value $2^k$, we have that $|X^c\cap X| = \infty$ for any conjugator $c$, since many of these generators will be taken to each other in a small neighbourhood of $0$ by the conjugation action of the first linear part of $c$.

Therefore $D$ is a normalish amenable subgroup of $F$.
\end{proof}

Similarly, we will now show that the standard version of R. Thompson's group $F$ in $T$ is normalish in $T$.
\begin{proof}
Let $K=\{k_1,k_2,\ldots k_m\}\subset T$ and consider the set 
\[
Y\seteq \{F^{k_i}\mid 1\leq i\leq k_m\}
\] which consists of subgroups of $T$ which are conjugates of $F$.

From the definition of $F$ given in the introduction, $F$ contains all orientation-preserving dyadic pl-homeomorphisms of $[0,1]$; these are the homeomorphisms which fix $0$ and $1$ (identified in the circle $S^1$), having (only finitely many) breaks in slope, all of which are occurring at dyadic rationals, and where all slopes of affine components are integral powers of $2$.  In particular, for any dyadic interval $(\frac{a}{2^j},\frac{b}{2^k})$, consider the set $B_{(\frac{a}{2^j},\frac{b}{2^k})}$ of all elements of $F$ supported exactly on $(\frac{a}{2^j},\frac{b}{2^k})$, which set is commonly known to generate a subgroup of $F$ isomorphic to $F$ (this follows easily from Cannon, Floyd and Parry's survey \cite{CFP}).  It is now immediately the case that the intersection
\[
\mathscr{I}\seteq \cap_{1\leq i\leq m}F^{k_i}
\]
is an infinite set.   This is easy to see, since it contains any such set $B_{(\frac{a}{2^j},\frac{b}{2^k})}$ where the interval $(\frac{a}{2^j},\frac{b}{2^k})$ is disjoint from the set $$\{0\cdot k_i\mid 1\leq i\leq m\}\cup\{0\}.$$
\end{proof}

We comment that there are other arguments for our Remark \ref{TNormalishIfFAmenable} which might be easier than the above, but, one should be careful:  $T$ admits subgroups with global fixed point set not empty, which are not conjugate in $T$ to subgroups of the canonical embedding of $F$ in $T$.

\end{document}